\documentclass[a4paper, 11 pt, twoside]{amsart}
\usepackage{srollens-en, shortcuts-sr}
\usepackage[left = 3.5cm, right = 3.5cm, headsep = 6mm,
footskip = 10mm, top = 35mm, bottom = 35mm, footnotesep=5mm, headheight =
2cm]{geometry}

\newlist{prooflist}{description}{1}
\setlist[prooflist]{font=\normalfont \itshape, labelindent = \parindent, leftmargin = 0pt}

\usepackage[unicode,bookmarks, pdftex]{hyperref}
\hypersetup{colorlinks=true,citecolor=NavyBlue,linkcolor=NavyBlue,urlcolor=Orange, pdfpagemode=UseNone, breaklinks=true}

\newcommand{\muhat}{\widehat{\mu}}
\DeclareMathOperator{\mult}{mult}
\DeclareMathOperator{\vol}{vol}

\newcommand{\Sig}{\overline{\Sigma}}
\newcommand{\Tot}{\operatorname{Tot}}

{\theoremstyle{Lehn-it}
\newtheorem{mainthm}{Theorem}
}

\title[From minimal valuations to irrational Seshadri constants]{Transforming local positivity: from minimal valuations to irrational Seshadri constants}

\author{S\"onke Rollenske}
\address{FB 12/Mathematik und Informatik\\
Philipps-Universit\"at Marburg\\Hans-Meerwein-Str. 6\\
35032 Marburg\\Germany}

\email{rollenske@mathematik.uni-marburg.de}

\begin{document}
\begin{abstract}
	We explain how a valuative Nagata type statement proposed by  Dumnicki, Harbourne, K\"uronya, Ro\'e and Szemberg on a polarised surface $(S, H)$ implies maximality (and hence often irrationality)  of Seshadri constants on simple cyclic covers. 
	
	As sample applications we show that very general degree two K3 surfaces and very general hypersurfaces in $\IP^3$ of degree $6$, $8$ or $d\geq 10$ have irrational Seshadri constant at a very general point, unless $d$ is a square. 
	
\end{abstract}


\maketitle
\setcounter{tocdepth}{1}
\tableofcontents

\section{Introduction}\label{sec:intro}
Let $X$ be a normal projective surface, $L$ an ample line bundle on $X$ and $x\in X$ a smooth point. The Seshadri constant
\[
\varepsilon(L;x)=\inf_{C\ni x}\frac{L\cdot C}{\mult_x C}
\]
measures the local positivity of $L$ at $x$. It satisfies $\varepsilon(L;x)\le\sqrt{L^2}$, and if the inequality is strict, then $\varepsilon(L;x)$ is computed by a so-called submaximal curve and is rational \cite[Prop.~2.1.1, Rem.~2.1.4]{BDHKKSS}. So intuitively a maximal Seshadri constant says that there are no unexpectedly singular curves through $x$. 
 Whether irrational Seshadri constants exist at all was a long-standing question and, in joint work with Alex K\"uronya, we could only get arbitrarily close \cite{KR}. The first examples were found only recently, by Laface and Ugaglia \cite{LU} and by Malara, Merta, Szpond and Zieli\'nski \cite{MMSZ}, on rational surfaces of Picard number at least ten.

At about the same time OpenAI released a proof of Nagata's conjecture \cite{OAI1}, which is about a different kind of local positivity on blow ups of the plane $\IP^2$. Tevelev's account of that proof \cite{Tev} also establishes the case of rational $t>9$ of a valuative version of Nagata's conjecture due to Dumnicki, Harbourne, K\"uronya, Ro\'e and Szemberg  \cite[Conj.~4.3]{DHKRS} (see also \cite{GMMM}):
for a curvilinear subscheme $Z\subset \IP^2$ of length $q$ let $\nu_Z$ be the monomial valuation of weights $(1,q)$ that measures contact with $Z$, see Section~\ref{sec:val} for details.
\begin{thm}[{\cite[Thm.~C]{Tev}, \cite{DHKRS} for $q=9$}]\label{thm: Tevelev plane DHKRS}
	Let $(S, H)= (\IP^2, \ko_{\IP^2}(1))$. Then for a very  general curvilinear subscheme of length $q\geq 9$ the valuation is minimal, that is,
	\begin{equation}\label{eq: minimal}
		\muhat_H(\nu_Z)=\sqrt{qH^2}.
	\end{equation}
\end{thm}
As suggested already in \cite{DHKRS}, this makes sense not only for the plane, so we pose the following definition.
\begin{defin}
	A normal polarised surface $(S,H)$ is \emph{DHKRS in degree $q$} (or just \emph{$q$-DHKRS}) if \eqref{eq: minimal} holds for a very general curvilinear subscheme $Z\subset S$ of length $q$.
\end{defin}
We note that $(S, H)$ is DHKRS in degree $1$ if and only if $\varepsilon(H;x)=\sqrt{H^2}$ at a very general point (Lemma~\ref{lem:width}). 

The main observation of this paper is, that such Nagata type statements can be transformed to maximality, thus often irrationality, of Seshadri constants via a cyclic cover construction.
\begin{mainthm}\label{thm:A}
	Let $(S,H)$ be a normal projective surface with an ample line bundle that is DHKRS in degree $q\ge2$. Let $N$ be a Cartier divisor on $S$ such that
	\begin{itemize}
		\item the general $\Delta \in |qN|$ is reduced,
		\item $|qN|$ is $q$-inclusive, that is, a general curvilinear $Z$ of length $q$ is contained in the smooth locus of a $\Delta \in |qN|$. 
	\end{itemize}	
		 Then for a very general $\Delta\in|qN|$ the cyclic cover $\pi\colon X\to S$ of degree $q$ branched along $\Delta$ satisfies
	\[
	\varepsilon(X,\pi^*H;x)=\sqrt{qH^2}\qquad\text{for very general }x\in X.
	\]
\end{mainthm}

For applications we need more DHKRS-surfaces and a minor modification of Tevelev's proof of Theorem \ref{thm: Tevelev plane DHKRS} provides a good starting list.
\begin{mainthm}\label{thm:B}
	The polarised surfaces $(S,H)$ in Table~\ref{tab: DHKRS} are DHKRS in the degrees $q$ listed there.
\end{mainthm}

\begin{table}
	\begin{tabular}{@{}clcl@{}}
		\toprule
		$S$ & $H$ & $H^2$ & DHKRS in degree\\
		\midrule
		$\IP^2$ & $\ko(1)$ & $1$ & $q\in\{1,4\}$ or $q\ge9$ \\
		$dP_\ell$, $1\le\ell\le9$ & $-K_S$ & $\ell$ & $q\ge\ell$\\
		$\IP^1\times\IP^1$ & $\ko(1,1)$ & $2$ & $q\ge8$\\
		$\IP(1,1,2)$ & $-K_S=\ko(4)$ & $8$ & $q\ge8$\\
		$\IP(1,2,3)$ & $-K_S=\ko(6)$ & $6$ & $q\ge6$\\
		$\Sym^2E$ & $\ko(C_p)$ & $1$ & $q\ge2$\\
		\bottomrule
	\end{tabular}
	\caption{Some known $q$-DHKRS surfaces}\label{tab: DHKRS}
	\end{table}

We only give two sample applications, using a simple propagation trick, deliberately leaving space for future explorations.
\begin{mainthm}\label{thm:C}
	In each of the following cases the Seshadri constant at very general points $x$ is maximal, $\varepsilon(L;x)=\sqrt{L^2}$:
	\begin{enumerate}
		\item very general polarised K3 surfaces $(X,L)$ of degree $L^2=2$;
		\item very general surfaces $X\subset\IP^3$ of degree $d\in\{4,6,8,9\}$ or $d\ge10$, with $L=\ko_X(1)$.
	\end{enumerate}
	In both cases, the very general such surface has Picard number $\rho (X) = 1$, and the Seshadri constant is irrational unless $d$ is a square in $(ii)$.
\end{mainthm}
To our knowledge, these are the first examples of surfaces with $\rho(X)=1$ and irrational Seshadri constants.

\subsection*{Acknowledgements}
I would like to thank Alex K\"uronya and Thomas Bauer for many discussions about Seshadri constants, in particular in the case of double covers of the plane, and Rita Pardini for encouraging comments on an earlier version.  I gratefully acknowledge support by the DFG.

\subsubsection*{Personal note}
Finding a surface with irrational Seshadri constant was a pet side interest of mine ever since I learned about the problem from Thomas Bauer. I pursued this in my work with Alex K\"uronya \cite{KR}, getting close to irrational but not really hitting home, and then further in an unsuccessful attempt to prove it for K3 surfaces of degree two with Alex and Thomas. After the OpenAI dump of the proof of the Nagata conjecture \cite{OAI1} I was curious if the problem was solved and learned about \cite{LU, MMSZ}. The first answer by AI \cite{Claude} to my question whether the K3 case was still open was an emphatic, \emph{yes, that is difficult}, but shortly afterwards it pointed out the calculation that became Proposition~\ref{thm:identity}. 

At that point, I regretted asking: it would have been so much more satisfying to find this myself. However, realistically speaking, the chances of me noticing it were minimal for two reasons: Seshadri constants are smaller at special points and I always assumed that points in the branch locus were too special in that sense (and for double covers of the plane they are). Also, knowing that there is a group of experts looking into the proof of the Nagata conjecture and being overwhelmed by the current flood of results anyway, studying the details of the proof and possible consequences was not high on my agenda. 

In the end, I thought the result was too interesting not to share and decided to write this note.

\subsection*{AI disclosure}
Claude \cite{Claude} pointed out the connection that leads to Theorem \ref{thm:A} and how to generalise Tevelev's work to get more DHKRS surfaces. Shaping of the material and deciding on the presentation was done by the author.

\section{Seshadri constants and monomial valuations}
Throughout this section $X$ is a normal projective surface, $L$ an ample line bundle on $X$ and $x\in X$ a smooth point. Recall, besides the facts mentioned at the beginning of the introduction, that $x\mapsto\varepsilon(L;x)$ is lower semicontinuous for the topology whose closed sets are the countable unions of Zariski closed sets, so that $\varepsilon(L;x)$ takes its maximal value at very general points \cite[2.2.8 and Sect.~7]{BDHKKSS}. 
This fact is usually stated for  smooth surfaces but either passing through a resolution or repeating the proof verbatim establishes it on a normal surface. 

\subsection*{Asymptotic maximal values}
Let $\nu$ be a real valuation of the function field of $X$ centred at a smooth point, and $D$ a line bundle. Following \cite[Sect.~2]{DHKRS} and \cite{BKMS} we put
\[
\mu_{kD}(\nu)=\max\{\nu(D'):D'\in|kD|\},\qquad \muhat_D(\nu)=\lim_{k\to\infty}\frac{\mu_{kD}(\nu)}k ,
\]
where $\nu(D')$ is the value of a local equation of $D'$. Since $\mu_{(k+l)D}(\nu)\ge\mu_{kD}(\nu)+\mu_{lD}(\nu)$, the limit exists and equals the supremum of the sequence. In particular
\begin{equation}\label{eq:sup}
	\nu(D')\le k\,\muhat_D(\nu)\qquad\text{for all }k\ge1\text{ and all }D'\in|kD|,
\end{equation}
and $\muhat_{cD}(\nu)=c\,\muhat_D(\nu)$ for $c\in\IZ_{>0}$. If $D$ is ample, then
\begin{equation}\label{eq:BKMS}
	\muhat_D(\nu)\ge\sqrt{D^2/\vol(\nu)},\qquad \vol(\nu)=\lim_{m\to\infty}\frac{\dim\ko_{X,x}/\{\nu\ge m\}}{m^2/2},
\end{equation}
by \cite[Prop.~2.9]{BKMS}; this just compares $h^0(mD)=m^2D^2/2+O(m)$ with the colengths of the valuation ideals. Following \cite{DHKRS}, we call $\nu$ \emph{minimal for $D$} if equality holds in \eqref{eq:BKMS}.

For $\nu=\ord_x$ we have $\vol(\ord_x)=1$, and $\muhat_L(\ord_x)$ is the asymptotic maximal order of vanishing at $x$ of sections of multiples of $L$. The following is probably well known to the experts.
\begin{lem}\label{lem:width}
	The following are equivalent:
	\begin{enumerate}
		\item $\varepsilon(L;x)=\sqrt{L^2}$;
		\item  $\ord_x$ is minimal for $L$, that is, $\muhat_L(\ord_x)=\sqrt{L^2}$;
		\item  $\ord_x\sigma\le k\sqrt{L^2}$ for all $k\ge1$ and all $0\ne\sigma\in H^0(X,kL)$.
	\end{enumerate}
\end{lem}

\begin{proof}
	$(ii)$ and $(iii)$ are equivalent by \eqref{eq:sup} and \eqref{eq:BKMS}.
	
	By the remark above we may assume that $X$ is smooth and $L$ nef and big. If $(iii)$ fails, some $\sigma\in H^0(X,kL)$ has $\ord_x\sigma>k\sqrt{L^2}$; writing $\divisor\sigma=\sum a_iC_i$ we get $\sum a_i\,L\cdot C_i=kL^2$, which is less than $\sqrt{L^2}\sum a_i\mult_xC_i$, so some $C_i\ni x$ has $L\cdot C_i<\sqrt{L^2}\mult_xC_i$, and $(i)$ fails.
	
	Conversely, suppose $(i)$ fails, so that there is an irreducible curve $C\ni x$ with $L\cdot C<\sqrt{L^2}\,r$, where $r=\mult_xC$. For small rational $a>0$ the $\IQ$-divisor $M=L-aC$ still satisfies $M^2>0$ and $L\cdot M>0$. Hence $\muhat_M(\ord_x)\ge\sqrt{M^2}$, by the count behind \eqref{eq:BKMS}: Riemann--Roch gives $h^0(kM)\ge k^2M^2/2+O(k)$, since $h^2(kM)=h^0(K_X-kM)=0$ for $k\gg0$ as $L$ is nef and eventually $L\cdot(K_X-kM) <0$. So for every $\tau<\sqrt{M^2}$ and $k$ sufficiently divisible there is a curve $C'\in|kM|$ with $\mult_xC'\ge k\tau$. Then $C'+kaC\in|kL|$ has multiplicity at least $k(\tau+ar)$ at $x$, and $(iii)$ fails as soon as $\sqrt{M^2}+ar>\sqrt{L^2}$. This holds for small $a$: as $M^2=L^2-2aL\cdot C+a^2C^2$, the derivative of $a\mapsto\sqrt{M^2}+ar$ at $a=0$ is $r-L\cdot C/\sqrt{L^2}>0$.
\end{proof}

To pass from special to very general surfaces we use the following semicontinuity; the special fibres may be singular, even non-normal, as long as the marked points are smooth.

\begin{lem}\label{lem:semicont}
	Let $f\colon\kx\to T$ be a flat projective morphism of irreducible varieties whose general fibre is irreducible, and $\kl$ an $f$-ample line bundle. Assume that the number of sections $h^0(X_t, \kl^{k}|_{X_t})$ on the fibre $X_t$ is independent of $t$ for every $k\ge1$. Let $\kx^\circ\subset\kx$ be the open set where $f$ is smooth. 
	If $\muhat_{\kl|_{X_{t_0}}}(\ord_{x_0})\le c$ for one point $x_0\in X_{t_0}\cap\kx^\circ$, then $\muhat_{\kl|_{X_t}}(\ord_x)\le c$ for very general $t\in T$ and very general $x\in X_t$.
	
	In particular, if $\ord_{x_0}$ is minimal and the general fibre is normal, then \[\varepsilon(\kl|_{X_t};x)=\sqrt{\kl|_{X_t}^2}\]
	 for  very general $t\in T$ and very general $x\in X_t$.
	
\end{lem}

\begin{proof}
	For $x\in X_t\cap\kx^\circ$ and $k\ge1$ put $\mu_k(x)=\mu_{\kl^k|_{X_t}}(\ord_x)$. Since $T$ is reduced and $h^0(X_t,\kl^k|_{X_t})$ is constant, $f_*\kl^k$ is locally free and $f_*\kl^k\otimes\kappa(t)\to H^0(X_t,\kl^k|_{X_t})$ is an isomorphism for every $t$ \cite[III.12.9]{Har}. Over $\kx^\circ$ the sheaf of relative $(r-1)$-jets of $\kl^k$ is locally free, and evaluation of sections is a map of vector bundles from $f^*f_*\kl^k$ to it, whose kernel at $x\in X_t$ consists of the sections of $\kl^k|_{X_t}$ vanishing to order at least $r$ at $x$. By semicontinuity of the rank, $\mu_k(x)\ge r$ is a closed condition on $\kx^\circ$, so
	\[
	U_k=\{x\in\kx^\circ : \mu_k(x)\le\mu_k(x_0)\}
	\]
	is open in $\kx$, and nonempty as $x_0\in U_k$. By \eqref{eq:sup} we have $\mu_k(x_0)\le kc$, hence $\muhat_{\kl|_{X_{f(x)}}}(\ord_x)=\sup_k\mu_k(x)/k\le c$ for all $x\in\bigcap_kU_k$.

	It remains to see that a very general $X_t$ meets each proper closed subset $W_k=\kx\setminus U_k$ in a proper closed subset. The set of $t$ with $X_t\subset W_k$ is $T\setminus f(\kx\setminus W_k)$, a proper closed subset of $T$, since flat morphisms of finite type are open.

	For the last claim note that $\kl|_{X_t}^2$ does not depend on $t$, as $f$ is flat; apply the first part with $c=\sqrt{\kl|_{X_t}^2}$, and then Lemma~\ref{lem:width}.
\end{proof}

\subsection{Curvilinear schemes and monomial valuations}\label{sec:val}
Let $S$ be a normal projective surface and $q\ge1$ an integer. We write $\kc_q(S)$ for the set of curvilinear subschemes $Z\subset S$ of length $q$ supported at a single smooth point $\bar x$. Such a $Z$ lies on a smooth curve germ $\{v=0\}$ at $\bar x$, and if $u$ restricts to a local parameter of this germ, then $Z=V(u^q,v)$. Thus $Z$ is determined by $\bar x$ and the $(q-1)$-jet of a smooth curve germ, and $\kc_q(S)$ is an irreducible variety of dimension $q+1$. Let $\nu_Z$ be the monomial valuation of weights $(1,q)$ in $(u,v)$,
\[
\nu_Z\Big(\sum_{a,b}c_{ab}\,u^bv^a\Big)=\min\{\,b+qa:\ c_{ab}\ne0\,\},
\]
which measures the contact of a germ with $Z$. If $\bar x$ is the origin of $\IA^2$, with coordinates $(x,y)$, and $v=y-\xi(x)$, $u=x$, then $\nu_Z$ is the quasimonomial valuation $v(\xi,q)$ with characteristic exponent $q$ of \cite[Def.~2.2]{DHKRS}.

\begin{lem}\label{lem:Z}
	\begin{enumerate}
		\item[(1)] For $n\ge0$ we have $\{F:\nu_Z(F)\ge qn\}=I_Z^n$. Hence $\nu_Z(F)=\max\{n:F^q\in I_Z^n\}$ depends only on $Z$.
		\item[(2)] The ideal $J_m(Z)=\{F:\nu_Z(F)\ge m\}$ has colength $\#\{(a,b)\in\IZ_{\ge0}^2:b+qa<m\}=m^2/2q+O(m)$. Thus $\vol(\nu_Z)=1/q$, and $\muhat_H(\nu_Z)\ge\sqrt{qH^2}$ for every ample line bundle $H$ on $S$.
	\end{enumerate}
\end{lem}

\begin{proof}
	(1) The ideal $\{\nu_Z\ge qn\}$ is generated by the monomials $u^bv^a$ with $b+qa\ge qn$, and each of them is divisible by $v^n$ if $a\ge n$, and by $u^{q(n-a)}v^a\in(u^q,v)^n$ if $a<n$. The reverse inclusion is clear. Since $\nu_Z$ is a valuation, $\nu_Z(F)\ge n$ if and only if $\nu_Z(F^q)\ge qn$, if and only if $F^q\in I_Z^n$. (2) The monomials $u^bv^a$ with $b+qa<m$ form a basis of the quotient; the last claim then follows from \eqref{eq:BKMS}.
\end{proof}

\begin{rem}[Very general valuations]\label{rem:verygeneral}
	Lemma~\ref{lem:semicont} and its proof hold verbatim with $\ord_x$ replaced by the valuation $\nu_Z$.
	In particular, $\muhat_H(\nu_Z)$ takes its smallest value for very general $Z\in\kc_q(S)$.
\end{rem}

We recall from the introduction the class of surfaces where through a very general curvilinear subscheme of length $q$ we have no curves of unexpected tangency.
\begin{defin}\label{def:DHKRS}
	Let $S$ be a normal projective surface, $H$ an ample line bundle on $S$ and $q\ge1$ an integer. We say that $(S,H)$ is \emph{DHKRS in degree $q$} (or just \emph{$q$-DHKRS}) if $\nu_Z$ is minimal for $H$ for very general $Z\in\kc_q(S)$, that is,
	\[
	\muhat_H(\nu_Z)=\sqrt{qH^2}.
	\]
\end{defin}

By \eqref{eq:sup}, $(S,H)$ is DHKRS in degree $q$ if and only if $\nu_Z(F)\le k\sqrt{qH^2}$ for very general $Z$, all $k\ge1$ and all $0\ne F\in H^0(S,kH)$. The condition only depends on the ray spanned by $H$, so it makes sense for ample $\IQ$-Cartier divisors as well. 

For $q=1$ we have $\nu_Z=\ord_{\bar x}$, and by Lemma~\ref{lem:width} the surface $(S,H)$ is DHKRS in degree one if and only if $\varepsilon(H;\bar x)=\sqrt{H^2}$ at very general points $\bar x$.

\section{Cyclic covers of DHKRS surfaces}

Let $S$ be a normal projective surface, $H$ an ample line bundle, $N$ a Cartier divisor on $S$, $q\ge2$, and $0\ne\sigma\in H^0(S,\ko_S(qN))$ with divisor $\Delta\in|qN|$. The \emph{(simple) cyclic cover of degree $q$ branched along $\Delta$} is the surface
\[
\pi\colon X=\{s^q=\sigma\}\subset\Tot(\ko_S(N))\to S
\]
in the total space of the line bundle $\ko_S(N)$, where $s$ is the tautological section of $\pi^*\ko_S(N)$ \cite[Sect.~I.17]{BHPV}. If $\Delta\ne0$ is reduced, then $X$ is a normal irreducible surface. The map $\pi$ is Galois with group $\mu_q$ acting by $s\mapsto\zeta s$, and the decomposition $\pi_*\ko_X=\bigoplus_{i=0}^{q-1}\ko_S(-iN)\,s^i$ into eigensheaves gives by the projection formula
\begin{equation}\label{eq:sections}
	H^0(X,\pi^*H^{k})=\bigoplus_{i=0}^{q-1}s^i\cdot \pi^*H^0\big(S,\ko_S(kH-iN)\big),
\end{equation}
where $\mu_q$ acts on the $i$-th summand by $\zeta^i$. In particular the invariant sections are the pull-backs of sections of $kH$.

Let $\bar x\in\Delta$ be a point at which both $S$ and $\Delta$ are smooth. If $v$ is a local equation of $\Delta$ and $(u,v)$ are local parameters at $\bar x$, then near $x_0=\pi^{-1}(\bar x)$ the surface $X$ is given by $s^q=v$. So $X$ is smooth at $x_0$, with local parameters $(u,s)$, and $\pi$ is $(u,s)\mapsto(u,s^q)$. The subscheme of length $q$ of $\Delta$ at $\bar x$ is $Z=V(u^q,v)\in\kc_q(S)$.

\begin{prop}\label{thm:identity}
	In this situation $\muhat_{\pi^*H}(\ord_{x_0})=\muhat_H(\nu_Z)$. Thus $\ord_{x_0}$ is minimal for $\pi^*H$ if and only if $\nu_Z$ is minimal for $H$; if $\Delta$ is reduced, this holds if and only if $\varepsilon(X,\pi^*H;x_0)=\sqrt{qH^2}$.
\end{prop}

\begin{proof}
	By \eqref{eq:sections} every $0\ne\tau\in H^0(X,\pi^*H^k)$ is a sum of eigensections $s^i\pi^*\tau_i$ with $\tau_i\in H^0(S,\ko_S(kH-iN))$. In the local analytic coordinates near $\bar x$ write $\tau_i=\sum_{a,b}c^{i}_{ab}u^bv^a$ for constants $c^{i}_{ab}$, so that $s^i\pi^*\tau_i=\sum_{a,b}c^{i}_{ab}u^bs^{qa+i}$. Different triples $(a,b,i)$ give different monomials, hence
	\begin{equation}\label{eq:pullback}
		\ord_{x_0}(\tau)=\min_{\tau_i\ne0}\ord_{x_0}(s^i\pi^*\tau_i),\qquad \ord_{x_0}(s^i\pi^*\tau_i)=i+\nu_Z(\tau_i).
	\end{equation}
	For $i=0$ this gives $\mu_{k\pi^*H}(\ord_{x_0})\ge\mu_{kH}(\nu_Z)$.
	
	Conversely, by the first equality in \eqref{eq:pullback} it suffices to bound the order of a single eigensection $s^i\pi^*\tau_i$. We cannot apply \eqref{eq:sup} to $\tau_i$ directly, since it is a section of $kH-iN$ rather than of a multiple of $H$. But the $q$-th power of an eigensection is invariant, hence a pull-back: as $s^q=\sigma$,
	\[
	(s^i\pi^*\tau_i)^q=\pi^*(\sigma^i\tau_i^q)\qquad\text{with}\qquad\sigma^i\tau_i^q\in H^0(S,qkH).
	\]
	By \eqref{eq:pullback} for pull-backs ($i=0$) and \eqref{eq:sup} we get
	\[
	q\,\ord_{x_0}(s^i\pi^*\tau_i)=\ord_{x_0}\pi^*(\sigma^i\tau_i^q)=\nu_Z(\sigma^i\tau_i^q)\le qk\,\muhat_H(\nu_Z).
	\]
	Hence $\mu_{k\pi^*H}(\ord_{x_0})\le k\,\muhat_H(\nu_Z)$ for all $k$, and the two invariants agree. By \eqref{eq:BKMS} and Lemma~\ref{lem:Z}, both have minimal value $\sqrt{qH^2}$. The last claim is Lemma~\ref{lem:width} on the now normal cover $X$.
\end{proof}

To apply Proposition~\ref{thm:identity} to a $q$-DHKRS surface we need branch divisors through a very general $Z\in\kc_q(S)$.
\begin{defin}\label{def:inclusive}
	Let $D$ be a Cartier divisor on $S$. The linear system $|D|$ is \emph{$q$-inclusive} if a general $Z\in\kc_q(S)$ is contained in a member of $|D|$ that is smooth at the support of $Z$.
\end{defin}
We provide an easy way to check this property.
\begin{lem}\label{lem:inclusive}
	Let $q\ge2$. For $Z\in\kc_q(S)$ supported at $\bar x$ let $Z'\supset Z$ be the subscheme of length $q+1$ with ideal $I_Z\cap\mathfrak m_{\bar x}^2$. If the restriction map $H^0(S,\ko_S(D))\to H^0(\ko_{Z'}(D))$ is surjective for some $Z\in\kc_q(S)$, then $|D|$ is $q$-inclusive.
\end{lem}

\begin{proof}
	The members of $|D|$ containing $Z$ form a linear subsystem of codimension $q$, those containing $Z'$ one of codimension $q+1$. A member containing $Z$ is singular at $\bar x$ if and only if its equation lies in $I_Z\cap\mathfrak m_{\bar x}^2$. So some member containing $Z$ is smooth at $\bar x$. Surjectivity is an open condition on $Z$, as the schemes $Z'$ form a flat family over the variety of jets of local parameters (with $I_{Z'}=(u^q,uv,v^2)$); so it holds for general $Z$.
\end{proof}

Now we are ready to prove the first main result.
\begin{proof}[Proof of Theorem \ref{thm:A}]
	Put $V=H^0(S,\ko_S(qN))$. The $Z\in\kc_q(S)$ for which $\nu_Z$ is minimal form a very general subset, and those that lie on a member of $|V|$ smooth at their support contain a dense open subset. So some $Z$, supported at $\bar x$, has both properties: $Z\subset\Delta_0=\divisor\sigma_0$ with $\Delta_0$ smooth at $\bar x$. Then $Z$ is the subscheme of length $q$ of $\Delta_0$ at $\bar x$, and by Proposition~\ref{thm:identity} the valuation $\ord_{x_0}$ is minimal for $\pi_{\sigma_0}^*H$, where $x_0=\pi_{\sigma_0}^{-1}(\bar x)$.
	
	The covers $X_\sigma$, $\sigma\in V\setminus0$, form a flat projective family over $V\setminus0$. Its total space is the cyclic cover of $S\times(V\setminus0)$ branched along the universal divisor, which is reduced as the general member of $|V|$ is; so the total space and the general $X_\sigma$ are normal and irreducible. By \eqref{eq:sections} the dimensions $h^0(X_\sigma,\pi_\sigma^*H^k)$ do not depend on $\sigma$. So Lemma~\ref{lem:semicont} applies and gives the claim.
\end{proof}

\section{Examples of DHKRS surfaces and inclusive linear systems}\label{sec:DHKRS}

This section presents the input from Tevelev's account \cite{Tev} of the proof of the Nagata conjecture, which provides us with the first examples of DHKRS surfaces.

\subsection{Warm-up: the plane via the valuative Nagata conjecture}

\begin{thm}\label{thm:plane}
	The plane $(\IP^2,\ko(1))$ is DHKRS in degree $q$ for $q\in\{1,4\}$ and for all $q\ge9$. It is not DHKRS in the degrees $q\in\{2,3,5,6,7,8\}$.
\end{thm}

\begin{proof}
	Recall from Section~\ref{sec:val} that $\nu_Z$ is the quasimonomial valuation $v(\xi,q)$ of \cite{DHKRS}, for a smooth curve germ $\xi\supset Z$.
	\begin{prooflist}
		\item[$q\in\{1,4,9\}$] This is \cite[Prop.~2.22]{DHKRS}; compare Lemma~\ref{lem:q=delta}.
		\item[$q\ge10$] By \cite[Thm.~C]{Tev}, with weights $(1,q)$, for very general jets of local parameters $(u,v)$ at a fixed point the monomial valuation $\nu$ of weights $(1,q)$ satisfies $\nu(F)<k\sqrt q$ for all $k\ge1$ and $0\ne F\in H^0(\IP^2,\ko(k))$. As $\PGL_3$ acts transitively on $\IP^2$,
		 this holds for $\nu_Z$ with $Z$ very general.
		\item[$q\in\{2,3,5,6,7,8\}$] For every $Z$ there is a curve $F$ of some degree $k$ with $\nu_Z(F)>k\sqrt q$: a line with $\nu_Z(F)\ge2$ for $q=2,3$, a conic with $\nu_Z(F)\ge5$ for $q=5,6$, a cubic with $\nu_Z(F)\ge8$ for $q=7$, and a sextic with $\nu_Z(F)\ge17$ for $q=8$. These exist because $\nu_Z(F)\ge m$ imposes $\#\{(a,b):b+qa<m\}<h^0(\ko(k))$ linear conditions (Lemma~\ref{lem:Z}).\qedhere
	\end{prooflist}
\end{proof}

\begin{lem}\label{lem:inclusiveP2}
	For $j\ge q\ge2$ the linear system $|\ko_{\IP^2}(j)|$ is $q$-inclusive.
\end{lem}

\begin{proof}
	Let $Z\in\kc_q(\IP^2)$ be supported at $\bar x$. In suitable affine coordinates centred at $\bar x$ we have $Z=V(x^q,\,y-\xi(x))$, where $\xi$ is a polynomial of degree less than $q$, the Taylor polynomial of a smooth curve germ containing $Z$. The curve $\{y=\xi(x)\}$ has degree at most $q-1<j$, contains $Z$ and is smooth at $\bar x$. Adding general lines, which do not pass through $\bar x$, we obtain a member of $|\ko(j)|$ that contains $Z$ and is smooth at $\bar x$.
\end{proof}

\subsection{Surfaces containing an elliptic curve}\label{sec:elliptic}
In this section we present a straightforward generalisation of \cite[Thm. C]{Tev} with nearly verbatim the same proof. 

For local parameters $(u,v)$ at a smooth point of a surface and positive integers $p,q,K$ let
\[
J^{p,q}_K(u,v)=\big(u^bv^a:\ pb+qa\ge K\big)
\]
be the ideal of the monomial valuation of weights $(p,q)$, extended by the unit ideal away from the point. Thus $J^{1,q}_K(u,v)=J_K(Z)$ for $Z=V(u^q,v)$, in the notation of Lemma~\ref{lem:Z}.

\begin{prop}[Tevelev]\label{prop:ruled}
	Let $\ell,p,q,h,K$ be positive integers with $q>\ell p$ and $K\ge h\sqrt{\ell pq}$. Then there is an elliptic curve $E$ with the following property: if $S$ is an integral projective surface and $C\cong E$ a curve in the smooth locus of $S$ with $C^2=\ell$, then there are a smooth point of $S$ and local parameters $(u,v)$ at it with
	\[
	H^0\big(S,\ko_S(hC)\otimes J^{p,q}_K(u,v)\big)=0.
	\]
\end{prop}

\begin{proof}
	For $S=\IP^2$ and $\ell=9$ this is \cite[Thm.~B and Prop.~4.1]{Tev}. We go through the three steps of the proof.
	\begin{prooflist}
		\item[The normal ruled surface] Let $M$ be a line bundle of degree $\ell$ on $E$, let $\Sig=\IP(\ko_E\oplus M)$ and $C_0\subset\Sig$ the zero section, so that $C_0^2=\ell$. Sections~2 and~3 of \cite{Tev} prove, for $\ell=9$ and a Tate curve $E$ with small parameter, that $H^0(\Sig,\ko(hC_0)\otimes J^{p,q}_K(w,y))=0$ for a suitable $M$, a point off the section at infinity and local coordinates $(w,y)$ there. The arguments hold verbatim with $9$ replaced by $\ell$ and $3h\sqrt{pq}$ by $h\sqrt{\ell pq}$: the results on theta functions \cite[Lemmas~3.1 and~3.2]{Tev} hold for line bundles of any degree, and in the rank estimate \cite[Prop.~2.4]{Tev} the hypotheses only enter through $\sigma=K/(ph)\ge\sqrt{\ell q/p}>\ell$, as in \cite[(2.10)]{Tev}.
		\item[The normal bundle] Translations act transitively on line bundles of degree $\ell$ on $E$, since $\theta_e^*M\cong M\otimes\ko_E(\ell(o-e))$ for the translation $\theta_e$ by $e$. They lift to isomorphisms of the ruled surfaces preserving $C_0$, so we may assume $M\cong\ko_S(C)|_C$.
		\item[Transfer to $S$] As in \cite[Prop.~4.1]{Tev}, degenerate $S$ to the normal cone of $C$ by blowing up $S\times\IA^1$ along $C\times\{0\}$. The proof there applies verbatim: it only uses that the exceptional divisor is $\Sig$ and the total space is smooth near it, which holds as $C$ lies in the smooth locus of $S$, and that sections of $\ko_S$ are constant.\qedhere
	\end{prooflist}
\end{proof}

\begin{thm}\label{thm:elliptic}
	Let $\ks\to B$ be a flat projective family of normal surfaces over an irreducible variety $B$, $\kh$ a line bundle on $\ks$ that is ample on the fibres, and $c\ge1$ an integer; put $\ell=c^2H_b^2$. Suppose that for every elliptic curve $E$ there are $b\in B$ and a smooth curve $C\in|cH_b|$ contained in the smooth locus of $S_b$ with $C\cong E$. Then for very general $b\in B$ the surface $(S_b,H_b)$ is DHKRS in degree $q$ for every $q>\ell$.
\end{thm}

\begin{proof}
	Fix $q>\ell$ and, for $h\ge1$, put $K_h=\left\lceil h\sqrt{\ell q}\right\rceil$. Let $\ku_h$ be the set of $(b,Z)$ in the irreducible relative variety $\kc_q(\ks/B)$ of curvilinear schemes with $H^0(S_b,H_b^{ch}\otimes J_{K_h}(Z))=0$. It is open by semicontinuity \cite[III.12.8]{Har}, as the sheaves $H_b^{ch}\otimes J_{K_h}(Z)$ form a flat family, and nonempty by Proposition~\ref{prop:ruled}, with $p=1$, applied to a curve $C\in|cH_b|$ isomorphic to the curve $E$ found there.

	For very general $(b,Z)$, hence for very general $b$ and very general $Z\in\kc_q(S_b)$, we have $(b,Z)\in\bigcap_h\ku_h$, that is, $\nu_Z(F)\le K_h-1<h\sqrt{\ell q}$ for all $h\ge1$ and $0\ne F\in H^0(S_b,H_b^{ch})$. So $\muhat_{cH_b}(\nu_Z)\le\sqrt{\ell q}=c\sqrt{qH_b^2}$, and equality holds by Lemma~\ref{lem:Z}.
\end{proof}
The degree $q=\ell$ can be added by an argument of \cite[Prop.~2.22]{DHKRS}.
\begin{lem}\label{lem:q=delta}
	Let $(S,H)$ be a normal polarised surface and $C\in|cH|$ an irreducible curve that is smooth at some smooth point of $S$. Then $(S,H)$ is DHKRS in degree $\ell=C^2$.
\end{lem}

\begin{proof}
By Remark~\ref{rem:verygeneral} and Lemma~\ref{lem:Z} it suffices to find one $Z\in\kc_\ell(S)$ with $\muhat_H(\nu_Z)\le\sqrt{\ell H^2}=cH^2$. Let $Z\subset C$ be the subscheme of length $\ell$ at a point where $C$ and $S$ are smooth, and choose local parameters with $C=\{v=0\}$, so that $\nu_Z(C)=\ell$. If a curve $F$ does not contain $C$, then the restriction of its local equation to $C$ vanishes to order at least $\nu_Z(F)$, so $\nu_Z(F)\le F\cdot C$. For $F\in|kH|$ write $F=mC+F'$ with $C\not\subset F'$; then
	\[
	\nu_Z(F)=m\ell+\nu_Z(F')\le m\ell+(k-mc)cH^2=kcH^2=k\sqrt{\ell H^2}.
	\]
	Thus $\muhat_H(\nu_Z)\le\sqrt{\ell H^2}$.
\end{proof}

\begin{proof}[Proof of Theorem \ref{thm:B}]
	The plane is Theorem~\ref{thm:plane}. Each of the other surfaces $(S,H)$ in Table~\ref{tab: DHKRS} belongs to a family $(S_b,H_b)_{b\in B}$ as in Theorem~\ref{thm:elliptic}: for every elliptic curve $E$ some $S_b$ contains a curve $C\cong E$ in its smooth locus with $C\in|cH_b|$. We give $B$, $c$ and $\ell=C^2$ below. Theorem~\ref{thm:elliptic} then gives the degrees $q>\ell$, and Lemma~\ref{lem:q=delta} the degree $q=\ell$.
	\begin{prooflist}
		\item[Del Pezzo surfaces] Let $B$ be the variety of $(9-\ell)$-tuples of points of $\IP^2$ in general position, and $c=1$. Every elliptic curve is a smooth plane cubic $C$, and general points $p_1,\dots,p_{9-\ell}$ of $C$ are in general position: by B\'ezout each of the conditions to be avoided (three points on a line, six on a conic, eight on a cubic singular at one of them) would force a linear equivalence on $C$, such as $p_i+p_j+p_k\sim\ko_C(1)$. The strict transform of $C$ lies in $|-K_S|$. Conversely, general points of $\IP^2$ lie on a smooth cubic, so a very general $S_b$ contains such a curve, as needed for Lemma~\ref{lem:q=delta}.
		\item[$\IP^1\times\IP^1$] Here $B$ is a point, $c=2$ and $\ell=8$. An elliptic curve $E$ has two non-isomorphic line bundles of degree two, and the two double covers $E\to\IP^1$ they define embed $E$ as a curve of bidegree $(2,2)$, since the image has arithmetic genus one.
		\item[$\IP(1,1,2)$] Here $B$ is a point, $c=1$ and $\ell=8$. Let $\IP(1,1,2)=\Proj\IC[x,y,w]$, so $-K=\ko(4)$. Given $E$, choose a line bundle $A$ of degree two, a basis $s_0,s_1$ of $H^0(E,A)$ and $s_2\in H^0(E,A^2)$ outside the span of $s_0^2,s_0s_1,s_1^2$. Then $(s_0:s_1:s_2)$ maps $E$ to $\IP(1,1,2)$, avoiding the vertex, and followed by the embedding by $\ko(2)$ it is the embedding by $|A^2|$. The image has degree two with respect to $\ko(1)$, and as $\ko(1)^2=\frac12$ it lies in $|\ko(4)|$.
		\item[$\IP(1,2,3)$] Here $B$ is a point, $c=1$ and $\ell=6$. On $\Proj\IC[z,x,y]$ with weights $1,2,3$ we have $-K=\ko(6)$, and every elliptic curve is a Weierstrass curve $y^2+a_1xyz+a_3yz^3=x^3+a_2x^2z^2+a_4xz^4+a_6z^6$, which avoids the two singular points.
		\item[$\Sym^2E$] Here $c=1$ and $\ell=1$. The Abel--Jacobi map $\Sym^2E\to\Pic^2(E)$ is a $\IP^1$-bundle with sections $C_p\cong E$ of self-intersection one, and $\ko(C_p)$ is ample \cite[V.2.21]{Har}. Let $B$ be a family of elliptic curves $E_b$ with a section $p_b$ in which every elliptic curve occurs, and $S_b=\Sym^2E_b$, $H_b=\ko(C_{p_b})$.\qedhere
	\end{prooflist}
\end{proof}

\subsection{Inclusive linear systems on del Pezzo surfaces}

\begin{lem}\label{lem:inclusivedP}
	Let $S$ be a del Pezzo surface of degree $\ell$ and $e\ge2$. Then $|-eK_S|$ is $e$-inclusive.
\end{lem}

\begin{proof}
	Let $C\in|-K_S|$ be smooth, $\bar x\in C$ a general point, and $Z\subset C$ the subscheme of length $e$ at $\bar x$. The restriction map $H^0(S,-eK_S)\to H^0(C,-eK_S|_C)$ is surjective, since $H^1(S,-(e-1)K_S)=0$ by Kodaira vanishing. On the elliptic curve $C$ the line bundle $M=-eK_S|_C$ has degree $e\ell$, so $h^0(C,M)=e\ell$, and by Riemann--Roch $h^0(C,M(-Z))=e\ell-e$ unless $M(-Z)$ is trivial. So $Z$ imposes $e$ independent conditions on $H^0(C,M)$, hence on $H^0(S,-eK_S)$. The exception $M\cong\ko_C(e\bar x)$ can only occur for $\ell=1$, and then only for finitely many $\bar x$. The member $C+D$ of $|-eK_S|$, with $D\in|-(e-1)K_S|$ not passing through $\bar x$, contains $Z$ and is smooth at $\bar x$. Hence the scheme $Z'$ of Lemma~\ref{lem:inclusive} imposes $e+1$ conditions, and the lemma applies.
\end{proof}

\section{Some surfaces with maximal Seshadri constant}\label{sec:appl}

Here is the first, maybe not so exciting example that follows directly from the solution to the valuative Nagata conjecture. 
\begin{exam}\label{exam:Xdm}
	Let $d=4$ or $d\geq 9$, and $m\geq 1$. Consider simple cyclic covers $X=X_{d,m}\to \IP^2$ of degree $d$ branched over a curve $\{f=0\}$ of degree $md$. We can realise them as hypersurfaces
	\[ X_{d,m} = \{ w^d - f(x,y,z) = 0 \} \subset \IP(1,1,1,m),\]
	and $L=\ko_X(1)$ is the pull-back of a line, so $L^2=d$.
	
	Since the plane is DHKRS in degree $d$ (Theorem~\ref{thm:plane}), and $|\ko(dm)|$ is $d$-inclusive (Lemma~\ref{lem:inclusiveP2}), with smooth general member, Theorem \ref{thm:A} applies and $\varepsilon(L;x)=\sqrt d$ at very general points $x$ of a very general $X_{d,m}$.
	
	This improves the bound $\sqrt d-d/m$ proved for $m\ge3$ in \cite{KR}, where it was also shown that the very general such  $X_{d,m}$ has Picard rank one.
\end{exam}

\begin{rem}[Propagation trick]\label{rem: propagation}
To get to the examples stated in Theorem \ref{thm:C} we use  a simple propagation trick that follows directly from Lemma \ref{lem:semicont}: 
if a single member of a family as in Lemma~\ref{lem:semicont} has a smooth point $x_0$ at which $\ord_{x_0}$ is minimal, for instance a point with maximal Seshadri constant, then the very general member of the family has maximal Seshadri constant at very general points. 

For us, the single member could be a very general cyclic cover or a birational image of one. 
\end{rem}

We  now prove Theorem \ref{thm:C} in the next two subsections.

\subsection{K3 surfaces of degree two}\label{sec:K3}
For double covers of the plane $X \to \IP^2$, the points on the ramification locus have lower Seshadri constants: the pullback of a tangent to the branch curve is a submaximal curve. So we need to find a different model, where Theorem \ref{thm:A} applies.

A del Pezzo surface of degree one is a sextic surface
\[
S=\{w^2=t^3+\alpha(x,y)\,t+\beta(x,y)\}\subset\IP(1,1,2,3),
\]
where $x,y,t,w$ have weights $1,1,2,3$ and $\alpha,\beta$ are  of degrees $4$ and $6$; moreover $-K_S=\ko_S(1)$ \cite[Ch.~8]{Dol}. The linear system $|-2K_S|$ is spanned by $x^2,xy,y^2,t$, so it is base point free and its general member is $\Delta=\{t=c(x,y)\}$ for a quadric $c$. The double cover $\pi\colon X\to S$ branched along $\Delta$, with $N=-K_S$, is obtained by adjoining a variable $z$ of weight one with $z^2=t-c$. Eliminating $t$, it is the double plane
\[
X=\{w^2=f(x,y,z)\}\subset\IP(1,1,1,3),\qquad f=(z^2+c)^3+\alpha\,(z^2+c)+\beta,
\]
with $\pi(x:y:z:w)=(x:y:z^2+c:w)$ and $\pi^*(-K_S)=\ko_X(1)$. So double planes whose branch sextic is even in $z$ are double covers of del Pezzo surfaces of degree one. Note that the map to the plane  $X\to\IP^2$ is \'etale at general points of the ramification curve $\{z=0\}$ of $\pi$.

\begin{proof}[Proof of Theorem~\ref{thm:C}$(i)$]
	Let $S$ be a very general del Pezzo surface of degree one. Then $(S,-K_S)$ is DHKRS in degree two (Theorem~\ref{thm:B}), and $|-2K_S|$ is $2$-inclusive (Lemma~\ref{lem:inclusivedP}) and base point free. By Theorem~\ref{thm:A}, for very general $\Delta\in|-2K_S|$ the double cover $X$ is smooth and $\varepsilon(\ko_X(1);x)=\sqrt2$ at very general $x\in X$, that is, $\ord_x$ is minimal (Lemma~\ref{lem:width}). By the above, $X$ is a double plane branched along a sextic. The double planes $\{w^2=f\}$, $\deg f=6$, form a family as in Remark~\ref{rem: propagation}, and its very general member is a very general polarised K3 surface of degree two, which concludes the proof.
\end{proof}

\subsection{Surfaces in \texorpdfstring{$\IP^3$}{P3}}\label{sec:P3}
The surfaces of degree $d$ in $\IP^3$ form a family as in Remark~\ref{rem: propagation}, so it suffices to find one surface of degree $d$ with a smooth point $y$ at which $\ord_y$ is minimal for $\ko(1)$. For $d\in\{4,9\}$ and $d\ge10$ the cyclic covers $X_{d,1}\subset\IP^3$ of Example~\ref{exam:Xdm} do the job. For $d=6,8$ we use birational images of cyclic covers of del Pezzo surfaces of degree two.

A del Pezzo surface of degree two is a quartic surface
\[
S=\{w^2=g(x,y,z)\}\subset\IP(1,1,1,2),
\]
with $-K_S=\ko_S(1)$; the anticanonical map is the double cover $S\to\IP^2$ branched along the quartic $\{g=0\}$, with covering involution $\iota\colon w\mapsto-w$. For $e\ge2$ we have $H^0(S,-eK_S)=\IC[x,y,z]_e\oplus w\cdot\IC[x,y,z]_{e-2}$.

\begin{lem}\label{lem:image}
	Let $e\ge2$ and $\sigma=a+wb\in H^0(S,-eK_S)$ with $a\in\IC[x,y,z]_e$ and $0\ne b\in\IC[x,y,z]_{e-2}$. Let $\pi\colon X\to S$ be the cyclic cover of degree $e$ branched along $\divisor\sigma$, with $N=-K_S$, and $L=\pi^*(-K_S)$. Then $H^0(X,L)$ has the basis $x,y,z,s$, where $s^e=a+wb$, and
	\[
	\varphi=(x:y:z:s)\colon X\to\IP^3
	\]
	is a finite birational morphism onto the surface $Y=\{(s^e-a)^2=b^2g\}$ of degree $2e$, with $\varphi^*\ko_Y(1)=L$.
\end{lem}

\begin{proof}
	By \eqref{eq:sections}, $H^0(X,L)=H^0(S,-K_S)\oplus s\cdot H^0(S,\ko_S)$. The sections $x,y,z$ have no common zero, so $\varphi$ is a morphism, and it is finite because $L$ is ample. On $X$ we have $wb=s^e-a$ and $w^2=g$, so $\varphi(X)\subset Y$. A general point $P\in\IP^2$ has two preimages $\tilde P$ and $\iota(\tilde P)$ in $S$, and the $2e$ points of $X$ over them have as $s$-coordinates the $e$-th roots of $(a+wb)(\tilde P)$ and of $(a-wb)(\tilde P)$. These two numbers are nonzero and different, as $w(\tilde P)b(P)\ne0$ for general $P$, so the $2e$ values of $s$ are distinct. Hence $\varphi$ is birational onto its image, which therefore has degree $L^2=2e=\deg Y$; so $\varphi(X)=Y$.
\end{proof}

\begin{proof}[Proof of Theorem~\ref{thm:C}$(ii)$]
	By the above, only $d=2e\in\{6,8\}$ remains. Let $S$ be a very general del Pezzo surface of degree two. Then $(S,-K_S)$ is DHKRS in degree $e$ (Theorem~\ref{thm:B}), and $|-eK_S|$ is $e$-inclusive (Lemma~\ref{lem:inclusivedP}) and base point free. By Theorem~\ref{thm:A} and Lemma~\ref{lem:width}, for very general $\sigma=a+wb$ the cyclic cover $X$ is smooth and $\ord_x$ is minimal for $L=\pi^*(-K_S)$ at very general $x\in X$. Let $\varphi\colon X\to Y$ be as in Lemma~\ref{lem:image}; note $b\ne0$. The surface $Y$ is not normal, as it is singular along the curve $\{b=0,\ s^e=a\}$ where the two sheets $w=\pm\sqrt g$ meet. But $\varphi$, being finite and birational, is an isomorphism over the normal locus of $Y$. So for very general $x$ the point $y=\varphi(x)$ is a smooth point of $Y$, and $\varphi^*\colon H^0(Y,\ko_Y(k))\to H^0(X,L^k)$ is injective and preserves the orders of vanishing at $y$ and $x$. Hence $\muhat_{\ko_Y(1)}(\ord_y)\le\muhat_L(\ord_x)=\sqrt d$, that is, $\ord_y$ is minimal for $\ko_Y(1)$.
\end{proof}

By the Noether--Lefschetz theorem \cite{Voi} very general surfaces of degree $d\ge4$ in $\IP^3$ have Picard number one. So we obtain irrational Seshadri constants $\sqrt d$ on surfaces of Picard number one for $d=6,8$ and for all non-square $d\ge10$.

\section{Remarks and open questions}\label{sec:questions}

\begin{rem}[Hypersurfaces in $\IP^3$ of low degree]
Theorem \ref{thm:C} does not cover hypersurfaces in $\IP^3$ of degree $1$, $2$, $3$, $5$, and $7$. 

For $d=1$ the maximality is trivial, while for $d=2$ the lines through $x$, and for $d=3$ the tangent plane section at $x$, which has a double point at $x$, show $\varepsilon(\ko_X(1);x)\le1$, respectively $\le\frac32$. 

The degrees $5$ and $7$ are still open, but one is led to believe that the Seshadri constant should be maximal for a very general point on a very general surface in these families. Our method does not reach them: if $L=\pi^*H$ for a cyclic cover $\pi\colon X\to S$ of degree $q$ and $(S,H)$ is DHKRS in degree $q$, then $L^2=qH^2$, so a prime $L^2$ forces $q=L^2$ and $H^2=1$. Among the surfaces of Theorem~\ref{thm:B} this leaves the plane, which is not DHKRS in degrees $5$ and $7$, and del Pezzo surfaces of degree one and $\Sym^2E$. In the latter case, a cover $(X, L)$ would have $h^0(X,L)\le3$ by \eqref{eq:sections}, so that $X$ does not map birationally to $\IP^3$.
\end{rem}

\begin{rem}[Cyclic covers of the plane of low degree]
	Example \ref{exam:Xdm} covers cyclic covers $X_{d,m}\to \IP^2$ of degree $d=4$ or $d\geq 9$ and Section \ref{sec:K3} covers the case $X_{2,3}$, from which the case $X_{2,m}$ for $m\geq 3$ can be deduced via Remark \ref{rem: propagation} by adding even multiples of a line to the branch locus. 
	
	Some other cases of low degree can probably be treated along the same lines as the hypersurfaces of degree $6$ in Section \ref{sec:P3}.
\end{rem}

\begin{rem}[Irregular surfaces]
	 If $E$ is a very general elliptic curve and $X\to\Sym^2E$ the cyclic cover of degree $q$ branched along a very general member of $|qnC_p|$, with $qn\ge3$, then $|qnC_p|$ is $q$-inclusive by the argument of Lemma~\ref{lem:inclusivedP}. Hence $\varepsilon(\pi^*C_p;x)=\sqrt q$ at very general points. These surfaces have irregularity one; they are properly elliptic for $(q,n)\in\{(2,2),(3,1)\}$ and of general type for $(q-1)n\ge3$.
\end{rem}

\begin{rem}[Variants of Theorem \ref{thm:A}]
	One could formulate more general versions of Theorem \ref{thm:A} for which the length of $Z$ does not coincide with the degree of the cover, or where the map is just finite,  but this is not needed for our application to Seshadri constants.
	
	Theorem \ref{thm:A} also holds with the same proof for a Weil divisor $N$. This could be useful for covers of the singular surfaces in Table~\ref{tab: DHKRS}.
	
	Similarly, the complete linear system $|qN|$ can be replaced by a linear subsystem $|V|\subset|qN|$ that is $q$-inclusive and whose general member is reduced.
	
\end{rem}

\begin{rem}[More DHKRS surfaces]\label{rem:moreDHKRS}
 Theorem~\ref{thm:elliptic} only applies to surfaces that are rational or ruled over an elliptic curve, since an elliptic curve $C\in|cH|$ in the smooth locus has $K_S\cdot C=-C^2<0$; among weighted projective planes only $\IP^2$, $\IP(1,1,2)$ and $\IP(1,2,3)$ contain such a curve, while other ruled surfaces over elliptic curves can be treated like $\Sym^2E$.
	
On the other hand, Claude suggests that the methods of \cite{OAI2} could possibly show that every surface $(S, H)$ is $q$-DHKRS for $q$ sufficiently large and it might be interesting to explore this further.
\end{rem}

\newcommand{\etalchar}[1]{$^{#1}$}

\end{document}